\documentclass[12pt]{amsart}

\usepackage{amssymb}

\usepackage{enumitem}

\makeatletter
\@namedef{subjclassname@2020}{%
  \textup{2020} Mathematics Subject Classification}
\makeatother

\newtheorem{theorem}{Theorem}[section]
\newtheorem{corollary}[theorem]{Corollary}
\newtheorem{lemma}[theorem]{Lemma}
\newtheorem{proposition}[theorem]{Proposition}

\theoremstyle{definition}
\newtheorem{definition}[theorem]{Definition}
\newtheorem{remark}[theorem]{Remark}

\numberwithin{equation}{section}

\def\sgn#1{\mbox{sgn}\left(#1\right)}
\def\norm#1{\|#1\|}

\def\nR{\mathbb R}
\def\nN{\mathbb N}
\def\beqa{\begin{eqnarray*}}
\def\eeqa{\end{eqnarray*}}
\def\beqal{\begin{eqnarray}}
\def\eeqal{\end{eqnarray}}
\def\nH{\mathcal H}
\def\nX{\mathcal X}
\def\nY{\mathcal Y}
\def\nZ{\mathcal Z}
\def\nA{A}
\def\tnorm#1{|\!|\!|#1|\!|\!|}
\def\dpX#1#2{\langle#1,#2\rangle_{\nX'\times \nX}}

\def\dpY#1#2{\langle#1,#2\rangle}
\def\normX#1{\norm{#1}_{\nX}}
\def\normY#1{\norm{#1}_{\nY}}
\def\dpXn#1#2{\langle#1,#2\rangle_{\nX_n'\times \nX_n}}
\def\scp#1#2{\langle#1,#2\rangle}
 
\def\normYS#1{\norm{#1}_{\nY'}}
\def\tnorm#1{\vert\!\vert\!\vert#1\vert\!\vert\!\vert}

\begin{document}


\baselineskip=17pt

\title[Mapped coercivity for nonlinear operators ]{The concept of mapped coercivity for nonlinear operators in Banach spaces}

\author{Roland Becker}
\address[R. Becker]{Department of Mathematics, Universit\'e de Pau et de l'Adour (UPPA), Avenue de l'Université, BP 1155, 64013 Pau CEDEX, France}
\email{becker@univ-pau.fr}

\author{Malte Braack}
\address[M. Braack]{Mathematical Seminar, Kiel University, Heinrich-Hecht-Platz 6, 24098 Kiel}
\email{braack@math.uni-kiel.de}

\date{\today}

\begin{abstract}
We provide a concise proof of existence for nonlinear operator equations in separable Banach spaces. Notably, the operator is not assumed to be monotone. Instead, our main hypotheses consist of a continuity assumption and a generalized coercivity property.
Mapped coercivity is a generalization of the usual coercivity property for nonlinear operators. In the case of linear operators,
we recover linear coercivity and the traditional inf-sup condition.
To illustrate the applicability of this general concept, we apply it to semi-linear elliptic problems and the Navier-Stokes equations.
\end{abstract}


\maketitle

%


\section{Introduction}\label{sec:1}

Let $\nX$ and $\nY$ be real, separable, and reflexive Banach spaces with dual spaces $\nX'$ and $\nY'$, respectively, and $\nA: \nX\to \nY'$. We are interested in the existence of solutions to the nonlinear operator equation
\beqal\label{eq:1}
   x\in\nX:\quad \nA(x) &=&  b\quad b\in \nY'.
\eeqal
In the classical Browder-Minty theory, see \cite{Browder1965, Browder1968}, the existence in the case $\nY=\nX$ relies on the coercivity property
\beqal\label{eq:standardcoercivity}
 \liminf_{x\in\nX, \normX{x}\to \infty}\frac{\dpY{\nA(x)}{x}}{\normX{x}} &=& +\infty,
\eeqal
and on monotonicity of the operator.
In this article we are interested in possible generalizations without monotonicity. Moreover, for indefinite operators $A$ the coercivity condition (\ref{eq:standardcoercivity}) cannot hold  and the following condition involving a mapping $\Phi:\nX\to\nY$ seems more appropriate:
\beqal\label{eq:mappedcoercivity}
 \liminf_{x\in\nX, \normX{x}\to \infty}\frac{\dpY{\nA(x)}{\Phi(x)}}{\normX{x}} &=& +\infty.
\eeqal
For linear operators $\nA$, such a condition has been proposed in \cite{BonnetBenChesnelCiarlet2012} with $\Phi= T\in L(\nX,\nY)$ (T-coercivity).
For linear $\Phi$, we can apply (\ref{eq:standardcoercivity}) to $\Phi'\circ\nA$.

The structure of this article is as follows: 
In Section~\ref{sec:2} we discuss several types of coercivity: without mappings in 
Section~\ref{sec:2-1}, with linear mappings in Section~\ref{sec:2-2}, and with nonlinear mappings in Section~\ref{sec:2-3}. The proofs are given in Section~\ref{sec:3}. 
In Section \ref{sec:4} we discuss the possibility to use duality maps as part of the
operator $\Phi$.
Finally, in Section \ref{sec:5}, we apply the theory to two partial differential equations, namely a class of semi-linear elliptic problems in the Sobolev spaces $W^{1,p}_0(\Omega)$ and the Navier-Stokes system in the Sobolev space $H^1(\Omega)$ for bounded Lipschitz domains $\Omega$.

\section{Mapped coercive operators}\label{sec:2}
We require weak continuity of $A$, i.e.
\beqal\label{eq:H1}
   x_n \rightharpoonup x\mbox{ in }\nX&\Longrightarrow&
   \nA(x_{n}) \stackrel{*}{\rightharpoonup} \nA(x)\mbox{ in }\nY'.
\eeqal
By the Uniform Boundedness Principle, the weak-$*$ convergence in (\ref{eq:H1}) implies that the images $(\nA(x_{n}))_{n\in\nN}$ are uniformly bounded in $\nY'$. 

In the case that $\nY$ is reflexive, the weak-$*$ topology coincides with the weak topology in $\nY$. 
Consequently, the weak-$*$-convergence in (\ref{eq:H1})  can be replaced by weak convergence. 
If, additionally, $A$ is compact and continuous w.r.t. the norm topologies of $\nX$ and $\nY'$, hypothesis  (\ref{eq:H1})  holds. 
However, it is important to consider that only weak-$*$-convergence, $\nA(x_{n}) \stackrel{*}{\rightharpoonup} \nA(x)$, is necessary for  (\ref{eq:H1}); the compactness of $A$ is not a requirement.

\subsection{Coercive operators}\label{sec:2-1}
In this section we state the results for coercive nonlinear operators. 

\begin{theorem}\label{lemma:1}
We assume that $\nA:\nX\to\nX'$ is  weakly continuous (\ref{eq:H1})
and coercive (\ref{eq:standardcoercivity}). Then,
 equation (\ref{eq:1}) has a solution.
\end{theorem}
The proof of this Theorem is given in Section \ref{sec:3-1} and \ref{sec:3-3}.

\subsection{Linear mappings}\label{sec:2-2}

In this section we state the results for nonlinear operators which are
coercive w.r.t. a linear mapping. 
The following property allows to obtain solutions for any data $b$.
\begin{definition}\label{def:1}
Let $\nX$ be reflexive.
We call $\nA:\nX\to\nY'$ {\em linearly mapped coercive}, 
if a surjection $\Phi\in L(\nX,\nY)$ exists 
s.t.
\beqal\label{eq:H2XX}
  \liminf_{x\in\nX, \normX{x}\to \infty}
  \frac{\dpY{\nA(x)}{\Phi x}}{\normX{x}} &=& +\infty.
\eeqal
\end{definition}
\begin{remark}
Clearly, $\Phi$ has to be a bijection.
In the case of a linear operator $A$, the definition of mapped coercivity (\ref{eq:H2XX}) is the so called {\em T-coercivity}
introduced by Bonnet-Ben, Chesnel and Ciarlet \cite{BonnetBenChesnelCiarlet2012,ChesnelCiarlet2013}.
\end{remark}
\begin{theorem}\label{thm:1}
We assume that $\nA:\nX\to\nY'$ is weakly continuous and linearly mapped coercive, 
i.e. (\ref{eq:H1}) and (\ref{eq:H2XX}) hold. Then, equation (\ref{eq:1}) has a solution.
\end{theorem}
\begin{proof}
Property (\ref{eq:H2XX}) implies that $\Phi^*A:\nX\to\nX'$ is coercive (\ref{eq:standardcoercivity}):
\beqa
  \liminf_{x\in\nX, \normX{x}\to \infty}
  \frac{\dpX{\Phi^*\nA(x)}{x}}{\normX{x}} &=& +\infty.
\eeqa
Since $\Phi^*b\in\nX'$, we obtain by Lemma~\ref{lemma:1} a solution $x\in\nX$ to 
\beqa 
     x\in\nX:\quad \Phi^*\nA(x) &=&  \Phi^*b.
\eeqa
Linearly mapped coercivity (\ref{eq:H2XX}) implies that $\Phi$ is bijective.
Hence, $\Phi^*\in L(\nY',\nX')$ is bijective. This implies that $x$ is solution to (\ref{eq:1}).
\end{proof}
The next definition allows for solutions corresponding to small data.
\begin{definition}\label{def:2}
We call the operator $A:\nX\to\nY$ {\em conditionally mapped coercive}, if there exists a surjection $\Phi\in L(\nX,\nY)$, s.t.
\beqal\label{eq:cond}
  M &:=& \liminf_{x\in\nX, \normX{x}\to \infty}\frac{\dpY{\nA(x)}{\Phi x}}{\normX{x}} > 0.
\eeqal
\end{definition}
\begin{corollary}\label{lemma:cond}
Let $\nA:\nX\to\nY'$ be weakly continuous (\ref{eq:H1}) and  conditionally mapped coercive (\ref{eq:cond}).
Let the right hand side be bounded by $\normYS{b} < M/\norm\Phi_{\nY\leftarrow\nX}$.
Then, equation (\ref{eq:1})  has a solution.
\end{corollary}
The proof of Corollary~\ref{lemma:cond} is given in Section \ref{sec:3-3}.

Let us now consider the particular case of a linear operator $A\in L(\nX,\nY')$. 
The next lemma states that, under further assumptions, the int-sup conditions, see Ne\v{c}as \cite{Necas1962} and Babu\v{s}ka \cite{Babuska1981},
 is equivalent to the linearly mapped coercivity assumption.
%
\begin{proposition}
Suppose that $\nY$ is a Hilbert space, and $A\in L(\nX,\nY')$. Then, the existence of $\gamma>0$ such that
\beqal\label{eq:infsup}
   \inf_{y\in\nY\setminus\{0\}} \sup_{x\in\nX\setminus\{0\}} \frac{\langle Ax,y\rangle}{\normX{x}\normY{y}} = \inf_{x\in\nX\setminus\{0\}} \sup_{y\in\nY\setminus\{0\}} \frac{\langle Ax,y\rangle}{\normX{x}\normY{y}} &=& \gamma,
\eeqal
is equivalent to (\ref{eq:H2XX}).
\end{proposition}
\begin{proof}
It is well known that (\ref{eq:infsup}) is equivalent to $A$ being an isomorphism.\\
$\Longrightarrow$:
Let $I_{\nY}\in L(\nY,\nY')$ be the  injective
antilinear operator isometry such that 
$\scp{y^*}{y}_{\nY'\times\nY} = (I_{\nY}^{-1}y^*,y)_{\nY}$.
If (\ref{eq:infsup}) holds, we define $\Phi := I_{\nY}^{-1} \circ A$, which is surjective. Then for any $x\ne0$
\[
\gamma\normX{x}\le\sup_{y\in\nY\setminus\{0\}}\frac{\scp{Ax}{y}}{\normY{y}} = \frac{\scp{Ax}{\Phi x}}{\normY{\Phi x}}\le \frac{\scp{Ax}{\Phi x}}{\normX{x}}\norm{\Phi^{-1}}_{\nX\leftarrow\nY}
\]
and (\ref{eq:H2XX}) is fulfilled. \\
$\Longleftarrow$:
Linearly mapped coercivity (Def. \ref{def:1})  implies 
\beqa
\alpha := \inf_{\norm{x}=1} \scp{Ax}{\Phi x} >0,
\eeqa
because, if there was $x\in\nX$ with $\normX{x}=1$ and $\scp{Ax}{\Phi x}<0$, then for $x_n=nx$ we have
\[
\frac{\scp{ \nA x_n}{\Phi x_n}}{\normX{x_n}} = n\scp{ \nA x}{\Phi x}\to-\infty.
\]
Moreover, if $(x_n)_{n\in\nN}\subset\nX$, $\normX{x_n}=1$ for all $n$, such that $\scp{ \nA x_n}{\Phi x_n}:=\lambda_n\to0$, 
then the sequence $w_n:=x_n/\lambda_n$ tends to infinity, but 
$\frac{\scp{ \nA w_n}{\Phi w_n}}{\norm{w_n}} = 1$. This contradicts Def. \ref{def:1} and thus $\alpha>0$.
It follows that $A^*\Phi$ is an isomorphism, so $\Phi$ is injective, and together with the assumption, is an isomorphism; so is then $A$.
But now, being $A$ an isomorphism implies (\ref{eq:infsup}).
\end{proof}

\subsection{Nonlinear mappings}\label{sec:2-3}

In the previous section, we saw that the presence of a linear map $\Phi:\nX\to\nY$ can be easily attributed to the 
coercive case by considering $\Phi^*A$. This is obviously not the case, if we want allow nonlinear $\Phi$.
In this section, we present a coercivity criterion with nonlinear mappings.
As before, the Banach space $\nX$ has to be  separable and reflexive. Furthermore, we assume the existence of a Schauder basis of $\nX$, and the dual $\nX'$ to be strictly convex. For any finite dimensional subspace $Z\subset\nX$, let $P_Z\in L(\nX,Z)$ be the natural projection induced by the Schauder basis.

The mapping $\Phi:\nX\to\nY$ we allow for is a composition of two
operators, $$\Phi=S\circ\Psi,$$ with $S\in L(\nX',\nY)$ and $\Psi:\nX\to\nX'$.
For $\Psi$ we need the following three properties to be satisfied:
\begin{enumerate}
\item Boundedness: It exists constant $C_\Psi>0$ s.t.
\beqal\label{eq:conditionX1}
  \norm{\Psi(x)}_{\nX'}  &\le&  C_\Psi  \normX{x} \quad\forall x\in\nX.
\eeqal
\item Definiteness: It holds
\beqal\label{eq:conditionX2}
    \dpX{\Psi(x)}{x}  >  0\quad\forall x\in\nX\setminus\{0\}.
\eeqal
\item Projection property: For any finite dimensional subspace $Z\subset\nX$ it holds
\beqal\label{eq:conditionX3}
    \dpX{\Psi(z)}{x-P_Zx}  &=&  0\quad\forall z\in Z,\ \forall x\in\nX.
\eeqal
\end{enumerate}
In Section \ref{sec:4} we will show that such an operator $\Psi:\nX\to\nX'$ exists. In particular, it can then be chosen as a scaled version of a duality map.
\begin{remark}
If $\Psi:\nX\to\nX'$ satisfies the properties (\ref{eq:conditionX1})-(\ref{eq:conditionX3}) and if $\sigma:\nX\to\nR_{>0}$
is a bounded function, $0<\sigma(x)\le M$, then any scaling $\widetilde\Psi:\nX\to\nX'$ defined by $\widetilde\Psi(x):=\sigma(x)\Psi(x)$
obviously also satisfies  (\ref{eq:conditionX1})-(\ref{eq:conditionX3}).
\end{remark}

\begin{definition}\label{def:3}
Let $\nX$ be reflexive.
We call $\nA:\nX\to\nY'$ {\em mapped coercive}, 
if there exists $\Psi:\nX\to\nX'$ with properties (\ref{eq:conditionX1})-(\ref{eq:conditionX3}) and a surjection $S\in L(\nX',\nY)$
s.t. for $\Phi:=S\circ\Psi$ it holds
\beqal\label{eq:H2XXX}
  \liminf_{x\in\nX, \normX{x}\to \infty}
  \frac{\dpY{\nA(x)}{\Phi (x)}}{\normX{x}} &=& +\infty.
\eeqal
\end{definition}
\begin{theorem}\label{thm:2}
We assume that $\nA:\nX\to\nY'$ is weakly continuous (\ref{eq:H1}) and mapped coercive (Def. \ref{def:3}). 
Then, equation (\ref{eq:1}) has a solution.
\end{theorem}
The proof is given in Section \ref{sec:3-3}.

\section{Proof of existence of solutions}\label{sec:3}
In this section we will give proofs of Theorem \ref{lemma:1} and Corollary \ref{lemma:cond}.
However, we start with proving existence of solutions for certain finite dimensional problems.

\subsection{Finite dimensional equations with coercive operators}\label{sec:3-1}
Since $\nX$ is separable, there exist a countable sequence of finite dimensional spaces $\nX_n:=\operatorname{span}\,(e_i)_{1\le i\le n}$ which is dense in $\nX$.
In this section, we prepare the proof of Proposition \ref{lemma:1} by considering
the discrete analog of (\ref{eq:1}), i.e a finite dimensional problem. 
For $\nY=\nX$, the discrete version reads:
\begin{equation}\label{eq:eqXnXn}
 x\in \nX_n:\quad \scp{A(x)-b}{y}_{\nX'\times\nX}=0\quad\forall y\in \nX_n,
\end{equation}
The restriction operator of functionals on $\nX$ to
functionals on $\nX_n$ will be denoted by $I_{\nX_n'}\in L(\nX,\nX_n')$.
\begin{lemma}\label{lemma:11}
Let $T_n\in L(\nX_n',\nX_n)$ be injective, and $R_n:\nX_n\to\nX_n$ given by
\beqal\label{eq:Rn}
  R_n(x)&:=& T_nI_{\nX_n'}(A(x)-b).
\eeqal
Then, $x\in\nX_n$ is solution of  (\ref{eq:eqXnXn})  iff $x\in N(R_n) $.
\end{lemma}
\begin{proof}
Since $T_n^*\in L(\nX_n',\nY_n)$ is surjective, we have the following equivalences:
\beqa
   x\in N(R_n) &\Longleftrightarrow&  T_nI_{\nX_n'}(A(x)-b)=0\\
    &\Longleftrightarrow&  \dpX{x'}{T_nI_{\nX_n'}(A(x)-b)} = 0\quad\forall x'\in\nX_n'\\
    &\Longleftrightarrow&  \dpY{I_{\nX_n'}(A(x)-b)}{T_n^*x'}_{\nY'\times\nY} = 0\quad\forall x'\in\nX_n'\\
    &\Longleftrightarrow&  \dpY{I_{\nX_n'}(A(x)-b)}{y}_{\nY'\times\nY} = 0\quad\forall y\in\nY_n\\
    &\Longleftrightarrow&  \dpY{A(x)-b}{y}_{\nY'\times\nY} = 0\quad\forall y\in\nY_n.
\eeqa
This completes the proof.
\end{proof}
Now, we consider the nonlinear eigenvalue problem of the operator defined in Lemma \ref{lemma:1}
with negative eigenvalue, i.e. the existence of $x\in\nX_n\setminus\{0\}$ and
$\lambda<0$ s.t.
\begin{equation}\label{eq:ev}
R_n(x) = \lambda x.
\end{equation}
Due to the finite dimension, $\nX_n$ is isomorph to its dual $\nX_n$.
Let $S_n\in L(\nX_n',\nX_n)$ be the corresponding isomorphism, i.e. for the basis $(e_i)_{i=1}^n$
of $\nX_n$ and for $x\in\nX_n$, $S_n^{-1}x$ is
the coordinate functional of $x$.
Due to the reflexivity of $\nX_n$, the adjoint can be considered as $S_n^*\in L(\nX_n',\nX_n)$
and is also an isomorphism. 
\begin{lemma}\label{lemma:2}
We assume that $\nA:\nX\to\nX$ is continuous on finite dimensional subspaces, and
coercive (\ref{eq:standardcoercivity}). Let 
 $R_n$ be given by (\ref{eq:Rn}) with $T_n:=S_n^*$.
Then, for any $n\in\nN$ there exists $r_n>0$ such that for any $\lambda<0$ and any solution $x\in\nX_n$ of 
(\ref{eq:ev}) with $R_n$ given by (\ref{eq:Rn}) it holds $\normX{x}< r_n$.
\end{lemma}
\begin{proof}
Since $S_n^{-*}\in L(\nX_n,\nX_n')$ is injective,, it holds for $x= \sum_{i=1}^n x_ie_i\not=0$:
\beqal\label{eq:conditionX3A}
   \scp{S_n^{-*} x}{x}_{\nX'\times \nX} & =&  \scp{S_n^{-1} x}{x}_{\nX'\times \nX} = \sum_{i=1}^n x_i^2> 0.
\eeqal
Suppose that $x\in\nX_n$ is a solution to (\ref{eq:ev})
with $\normX{x}>0$. 
With deduce
\beqa
0 &>&  \lambda\scp{S_n^{-*} x}{x}_{\nX'\times \nX}.
\\ &=& \scp{S_n^{-*}R_n(x)}{x}_{\nX'\times \nX}.
\\ &=& \scp{I_{\nX_n'}(A(x)-b)}{x}_{\nX'\times \nX}.
\\ &=& \scp{A(x)-b}{x}_{\nX'\times \nX}
\\ &\ge& \scp{A(x)}{x}_{\nX'\times \nX} - \norm{b}_{\nY'}\normX{x}.
\eeqa
However, for $r_n>0$ sufficiently large (independent of $\lambda$)
and $\normX{x}\ge r_n$, (\ref{eq:standardcoercivity}) implies that
the right hand side of the last estimate is positive, which is a contradiction. Hence, we derived $\normX{x}<r_n$.
\end{proof}

\begin{proposition}\label{prop:2}
We assume that $\nA:\nX\to\nX'$ is coercive (\ref{eq:standardcoercivity}) 
and  continuous on finite dimensional subspaces. Then,
 equation (\ref{eq:eqXnXn}) has at least one solution.
\end{proposition}
\begin{proof} 
(a) We choose $r_n>0$ as the radius in Lemma \ref{lemma:2}.
Suppose (\ref{eq:eqXnXn}) has no solution in $B_r:=\{x\in\nX_n: \normX{x}\le r_n\}$. Then $B_r\cap N(R_n)=\emptyset$ and the operator $T:B_r\to B_r$ with
\[
T(x) := - \frac{r_n}{\normX{R_n(x)}}R_n(x),
\]
 is well-defined. 
Since the restricted operator $A\vert_{\nX_n}$ was assumed to be continuous, $T$ is continuous.
Therefore, by Brouwer's fixed point theorem there exists a fixed point, $T(\bar x)=\bar x\in B_r$. Then, $\bar x$ is a solution to the eigenvalue problem (\ref{eq:ev}) with $\lambda=-\normX{R_n(\bar x)}/r_n<0$. Because of $\normX{\bar x}=r_n$, this contradicts Lemma \ref{lemma:2}.\\
(b) For the boundedness, we
suppose existence of an unbounded sequence $(x_n)_{n\in\mathbb N}$ of solutions.
By  equation (\ref{eq:eqXnXn}) we obtain the bound:
\beqa
  \frac{\dpY{\nA(x_n)}{x_n}}{\normX{x_n}} &=& 
  \frac{\dpY{b}{x_n}}{\normX{x_n}} \le 
  \norm{b}_{\nX'}
   \ <\ +\infty.
\eeqa
This gives for  the limes $n\to\infty$ a contradiction to coercivity  (\ref{eq:standardcoercivity}). 
\end{proof}

\subsection{Finite dimensional equations with mapped coercive operators}\label{sec:3-2}

In the case of a Schauder basis $(e_i)_{i\in\nN}$  of $\nX$ and  $\nX_n:=\operatorname{span}(e_i)_{1\le i\le n}$ let  $P_{n}:\nX\to \nX_n$ be the natural projections. 
The analogue of Lemma \ref{lemma:2} reads:
\begin{lemma}\label{lemma:ev}
We assume that $\nA:\nX\to\nY$ is continuous on finite dimensional subspaces, and 
mapped coercive (Def. \ref{def:3}). 
Let $R_n:\nX_n\to\nX_n$ be given by
\beqal\label{eq:Rn2}
  R_n(x)&:=& P_nS^* (A(x)-b).
\eeqal
Then, for any $n\in\nN$ there exists $r_n>0$ such that for any $\lambda<0$ and any solution $x\in\nX_n$ of 
(\ref{eq:ev}) it holds $\normX{x}< r_n$.
\end{lemma}
\begin{proof}
Suppose that $x\in\nX_n$ is a solution to (\ref{eq:ev})
with $\normX{x}>0$. 
With the coercivity property (\ref{eq:conditionX2})
and projection property (\ref{eq:conditionX3}), 
of $\Psi$ 
we deduce
\beqa
0 &>& \lambda \scp{\Psi(x)}{x}_{\nX'\times \nX}
\\ &=& \scp{\Psi(x)}{R_n(x)}_{\nX'\times \nX}
\\ &=& \scp{\Psi(x)}{P_nS^*(A(x)-b)}_{\nX'\times \nX}
\\ &=& \scp{\Psi(x)}{S^*(A(x)-b)}_{\nX'\times \nX}
\\ &=& \scp{A(x)-b}{S\Psi(x)}_{\nY'\times \nY}
\\ &\ge& \scp{A(x)}{S \Psi(x)}_{\nY'\times \nY} - c\normX{x},
\eeqa
with  $c>0$, depending on $\norm{b}_{\nY'}$, $C_S$ and the constants $C_\Psi$
in (\ref{eq:conditionX1}). However, for $r_n>0$ sufficiently large (independent of $\lambda$)
and $\normX{x}\ge r_n$, 
the right hand side of the last estimate
is positive, which is a contradiction. Hence, we derived $\normX{x}<r_n$.
\end{proof}
The finite dimensional spaces $\nY_n$ we are considering now are given by
\beqa
    \nY_n &:=& S\nX_n'.   
\eeqa
The following lemma ensures that this setting leads to dense subspaces of $\nY$.
\begin{lemma}\label{lemma:Xprimedense}
Let $S\in L(\nX',\nY)$ be surjective.
Then the images of the $S\nX_n'$ are dense in $\nY$, i.e.
\beqa
  \overline{\bigcup_{n\in\nN} S\nX_n'} &=& \nY.
\eeqa
\end{lemma}
\begin{proof}
Since $\bigcup_{n\in\nN} \nX_n$ is dense in $\nX$ and $\nX$ is reflexive, it follows that
$\bigcup_{n\in\nN} \nX_n'$ is dense in $\nX'$. This implies with the 
continuity and surjectivity of $S$:
\beqa
  \overline{\bigcup_{n\in\nN} S\nX_n'} &=& \overline{S\left(\bigcup_{n\in\nN} \nX_n'\right)} \ =\ 
  S\left(\overline{\bigcup_{n\in\nN} \nX_n'}\right) \ =\ S\nX'\ =\ \nY.
\eeqa
\end{proof}
\begin{proposition}\label{prop:dualmapped}
We assume that $\nA:\nX\to\nY$ is continuous on finite dimensional subspaces, and 
mapped coercive (\ref{eq:H2XXX}). Then, equation 
 \begin{equation}\label{eq:eqXnYn}
 x\in \nX_n:\quad \scp{A(x)-b}{y}_{\nX'\times\nX}=0\quad\forall y\in \nY_n,
\end{equation}
has a solution.
Furthermore, any sequence $(x_n)$ of solutions to (\ref{eq:eqXnYn}) is bounded.
\end{proposition}
\begin{proof}
(a) The proof for existence of solutions is identical to the proof of Proposition \ref{prop:2}, because we can use
Lemma \ref{lemma:ev}.\\
(b) For the boundedness of solutions we assume $x_n\in\nX_n$ to be a solution of (\ref{eq:eqXnYn}).
We have $P_n^*\Psi (x_n)\in\nX_n'$, and therefore  $SP_n^*\Psi (x_n)\in\nY_n$.
Now, we make two times use of the projection property (\ref{eq:conditionX3}):
\beqa
   \dpY{\nA(x_n)}{S\Psi (x_n)} &=& \dpX{\Psi (x_n)} {S^*\nA(x_n)}\\
   &=& \dpX{\Psi (x_n)} {P_nS^*\nA(x_n)}\\
   &=& \dpY{A(x_n)}{SP_n^*\Psi (x_n)} \\
   &=& \dpY{b}{SP_n^*\Psi (x_n)} \\
   &=& \dpX{\Psi (x_n)}{P_nS^*b} \\
   &=& \dpX{\Psi (x_n)}{S^*b} \\
   &=&  \dpY{b}{S\Psi (x_n)} .
\eeqa
Due to boundedness of $\Psi$ (\ref{eq:conditionX1}), we now have:
\beqa
  \frac{\dpY{\nA(x_n)}{S\Psi (x_n)}}{\normX{x_n}} &=& 
  \frac{\dpY{b}{S\Psi( x_n)}}{\normX{x_n}} \le 
  \normYS{b}\norm{S}_{\nY\leftarrow\nX}C_\Psi
   \ <\ +\infty.
\eeqa
The mapped coercivity (\ref{eq:H2XXX})  implies boundedness of the sequence $(x_n)$.
\end{proof}

\subsection{Proof of Theorems \ref{lemma:1}, \ref{thm:2} and Corollary \ref{lemma:cond}}\label{sec:3-3}


\begin{proof} (of Theorem \ref{lemma:1} and  \ref{thm:2})
The weak continuity (\ref{eq:H1}) imply that the finite dimensional restrictions $A\vert_{\nX_n}$ are continuous.
Both, coercivity and mapped coercivity, imply with Proposition \ref{prop:2} and Proposition \ref{prop:dualmapped}, respectively,
the existence of a bounded sequence of discrete solution $(x_n)_{n\in\mathbb N}$.
By the  Eberlein--\v{S}mulian theorem, see Diestel \cite{Diestel1984}, 
there exists a 
sub-sequence (denoted also by $(x_n)$) and an $x\in\nX$, such that $x_{n}\rightharpoonup x$. 
We know that   
$\bigcup_{n\in\mathbb N}\nY_n$ is dense in $\nY$: For coercive operators by separability of $\nY$, 
and for mapped coercive operators due to Lemma \ref{lemma:Xprimedense}.
So we can find for arbitrary $y\in \nY$
a sequence $(y_n)_{n\in\nN}$ with $y_n\in\nY_n$ and $y_n\to y$. 
We write
\beqa
\dpY{\nA(x)-b}{y} &=& \underbrace{\dpY{\nA(x_n)-b}{ y_n}}_{=0} + \underbrace{\dpY{\nA(x_n)}{ y-y_n}}_{\rm I} \\
&&+ \underbrace{\dpY{\nA(x)-\nA(x_n)}{y}}_{\rm II}  - \underbrace{\scp{b}{y-y_n}}_{\rm III}.
\eeqa
By continuity we have ${\rm III}=\scp{b}{y-y_n}\to 0$ for $n\to\infty$.
With  (\ref{eq:H1}) and the boundedness of $(\nA(x_{n}))_n$ it follows
\beqa
 \vert{ \rm I}\vert = \vert \dpY{\nA(x_{n})}{ y-y_{n}} \vert&\le& \normYS{\nA(x_{n})} \norm{y-y_{n}}_{\nY}\to0.
\eeqa
By (\ref{eq:H1}) it holds $\nA (x_{n}) \stackrel{*}{\rightharpoonup}\nA(x)$ in $\nY'$ which implies 
\beqa
{\rm II} =   \dpY{\nA(x_{n})-\nA(x)}{ y} &\to& 0.
 \eeqa
 This finishes the proofs of Theorems \ref{lemma:1} and \ref{thm:2}.
\end{proof}
\begin{proof} (of Corollary \ref{lemma:cond})
(a) We will firstly show the same result as in Lemma \ref{lemma:ev} under the
somehow weaker condition (\ref{eq:cond}). Let $n\in\nN$ be fixed. 
By assumption (\ref{eq:cond}), $\epsilon:=M- \normYS{b} \norm{\Phi_n}_{\nY\leftarrow \nX'}>0$.
As before, we have for any solution $x\in\nX_n$ of (\ref{eq:ev}):
\beqa
0 &>&  \scp{A(x)}{\Phi x}_{\nY'\times \nY} - \norm{b}_{\nY'}\norm{\Phi}_{\nY\leftarrow\nX}\normX{x}.
\eeqa
Due to (\ref{eq:cond}) there is  $r>0$ s.t. for any $x\in\nX_n$ with
$ \normX{x}\ge r$ it holds
\beqa
 \frac{\dpY{\nA(x)}{\Phi x}}{\normX{x}} &>& M-\epsilon = \normYS{b} \norm{\Phi}_{\nY\leftarrow \nX'}.
\eeqa
These two inequalities imply a contradiction. Hence, there exists a solution  $x\in\nX_n$ of (\ref{eq:ev}) and this is bounded
as $ \normX{x}< r$.\\
(b) Existence of solutions to (\ref{eq:eqXnYn}) and boundedness for $n\to\infty$ are obtained as in the proof of Proposition \ref{prop:dualmapped}.\\
(c) The existence of solutions to (\ref{eq:1}) now follows as in the proof above (of Theorem \ref{lemma:1} and \ref{thm:2}).
\end{proof}
%

\section{Duality maps}\label{sec:4}
In this section we will show that operators required in Def. \ref{def:3}, i.e. $\Psi:\nX\to\nX'$ with the properties (\ref{eq:conditionX1})-(\ref{eq:conditionX3}), exist, if $\nX'$ is strictly convex and $\nX$ has a Schauder basis.

\subsection{Duality maps w.r.t. equivalent norms}\label{sec:4-1}
Let $\tnorm\cdot$ an arbitrary norm on $\nX$ equivalent to $\normX\cdot$, and $\tnorm\cdot_{\nX'}$ its dual norm. 
The duality map $\Psi:\nX\to 2^{\nX'}$, normalized with respect to $\tnorm\cdot$, is
defined by
\beqal\label{eq:defJ}
  \Psi(x) &:=& \operatorname{argmax}\{\dpX{x'}{x}: x'\in\nX', \tnorm{x'}_{\nX'}\le\tnorm{x}\}.
\eeqal
In the following we give some well-known facts, see \cite{Cioranescu1990}. First, we apply for $x\in\nX$ the Hahn-Banach extension theorem to the functional $f_x$ on the 1-dimensional space $[x]:=span\langle x\rangle$:
\beqal\label{eq:toto11}
f_x :[x] \to \nR,\quad f_x(\lambda x):=\lambda\norm{x}^2,
\ \forall\lambda\in\nR,
\eeqal
which yields a maximizer in (\ref{eq:defJ}), since $\tnorm{f_x}_{[x]'}=\tnorm{x}$. It follows that the set $\Psi(x)$ is never empty. From this, we also see that the maximum 
in (\ref{eq:defJ}) is attained with $\dpX{x'}{x}=\tnorm{x'}_{\nX}\tnorm{x}=\tnorm{x}^2$, we obtain the equivalent characterization
\beqal\label{eq:defJ2}
  \Psi(x) = \{x'\in\nX': \tnorm{x'}_{\nX'}=\tnorm{x},\ \dpX{x'}{x}=\tnorm{x}^2 \}.
\eeqal
Multiplying the equations defining the set in (\ref{eq:defJ2}) by $\vert\lambda\vert$ and $\lambda^2$, respectively,  we get homogeneity of $\Psi$:
\beqa
\Psi(\lambda x)=\lambda \Psi(x)\quad\forall \lambda\in\nR.
\eeqa
This in turn shows that $\Psi(x)$ is the set of all norm-preserving extensions of $f_x$ defined in (\ref{eq:toto11}). Indeed, if $x'\in \Psi(x)$, we have
$
\dpX{x'}{x}=\tnorm{x}^2 = f_x(x)
$
and by homogeneity it follows that 
$\dpX{x'}{\lambda x} = \dpXn{\lambda x'}{ x} = f_x(\lambda x)$ and $x'$ is an extension of $f_x$.
Since the Hahn-Banach extension is unique on spaces with strictly convex dual space, we find for our
space $\nX_n$ that $\Psi(x)$ is single-valued.

We have the following properties of the duality map:
\begin{lemma}\label{lemma:X}
Let $\nX$ be a separable, reflexive Banach space with strictly convex dual $\nX'$.
Then,  the duality map $\Psi:\nX\to\nX'$ w.r.t to any norm $\tnorm\cdot$ equivalent to $\normX\cdot$ has the properties (\ref{eq:conditionX1}) and (\ref{eq:conditionX2}).
\end{lemma}
\begin{proof}
Due to norm equivalence and norm conservation of $\Phi$ we have for $x\in\nX$ and some constants $c_1,c_2>0$:
\beqal\nonumber
   c_1 \norm{\Psi(x)}_{\nX'} &\le&
    \tnorm{\Psi(x)}_{\nX'}  \ =\   \tnorm{x}  \le c_2  \normX{x} ,\\
    \dpX{\Psi(x)}{x}  &=&  \tnorm{x}^2\ge  c_1\normX{x}^2.\label{eq:psicoercive}
\eeqal
This shows the boundedness and definiteness.
\end{proof}
\begin{remark}
In the case that $\nX$ is a Hilbert space, with inner product $(\cdot,\cdot)_\nX$, the duality map $\Psi\in L(\nX,\nX')$ 
w.r.t. the norm $\normX\cdot$, are the injective
antilinear operator isometry:
\beqa
 \dpX{\Psi(x)}{x\rangle} &=&  (x,x)_\nX \ =\ \normX{x}^2\quad\forall x\in\nX.
\eeqa
\end{remark}

\subsection{Projection properties of duality maps}\label{sec:4-2}

\begin{lemma}\label{lemma:mapprojection}
Let $\nX$ be a Banach space with strictly convex dual $\nX'$, and $\tnorm\cdot$ a norm on $\nX$, equivalent to $\normX\cdot$. Furthermore, let for any $n\in\nN$ be $P_n:\nX\to\nX_n$ a projection with $\tnorm{P_n}_{\nX_n\leftarrow\nX}\le1$.Then
the duality map $\Psi:\nX\to\nX'$ w.r.t. $\tnorm\cdot$ satisfies the properties (\ref{eq:conditionX1})-(\ref{eq:conditionX3}).
\end{lemma}
\begin{proof}
According to Lemma \ref{lemma:X}, (\ref{eq:conditionX1})-(\ref{eq:conditionX2}) are valid for any duality map w.r.t. a norm equivalent to $\normX\cdot$. \\
(a) For verifying 
the projection property (\ref{eq:conditionX3}) we consider for $x'\in\nX_n'$ the following extension $E(x')\in\nX'$:
\beqal\label{eq:xPn}
    \dpX{E(x')}{x} &=& \dpX{x'}{P_nx}.
\eeqal
We have on the one hand:
\beqa
   \tnorm{E(x')}_{\nX'} &=& \sup_{x\in\nX} \frac{\dpX{x'}{P_nx}}{\tnorm{x}_\nX}\ \le\ \tnorm{x'}_{\nX'_n}\sup_{x\in\nX} \frac{\tnorm{P_nx}_\nX}{\tnorm{x}_\nX}.
\eeqa
Due to the assumption $\tnorm{P_n}_{\nX_n\leftarrow\nX}\le1$ we deduce
$\tnorm{E(x')}_{\nX'}\le\tnorm{x'}_{\nX'_n}$. Since $E(x')$ is an extension, it holds $\tnorm{E(x')}_{\nX'}\ge\tnorm{x'}_{\nX'_n}$.
This implies that $E(x')$ conserves the norm $\tnorm{E(x')}_{\nX'_n} =\tnorm{x'}_{\nX_n}$.\\
(b) Let now $z\in\nX_n$. Application of the  duality map $\Psi(z)$ leads to the unique norm conservative extension of the functional $f_z$
given as in (\ref{eq:toto11}), i.e. $\Psi(z)=E(f_z)$.
This implies for any $x\in\nX$:
\beqa
    \dpX{\Psi(z)}{x-P_nx}  &=&  \dpX{f_z}{P_n(x-P_nx)}\ =\ 0.
\eeqa
\end{proof}
One possibility to ensure $\tnorm{P_n}_{\nX_n\leftarrow\nX}\le1$ is the existence of a Schauder basis $(e_i)_{i\in\nN}$ for $\nX$.
Let $\nX_n:=\operatorname{span}(e_i)_{1\le i\le n}$, and let $P_{n}:\nX\to \nX_n$ be
 the natural induced projections.

\begin{corollary}\label{lemma:dualitymap}
Let $\nX$ be a Banach space with Schauder basis $(e_i)_{i\in\nN}$ and strictly convex dual $\nX'$. Then
the duality map $\Psi:\nX\to\nX'$ w.r.t. $\tnorm\cdot_\nX$, given by \beqal\label{eq:tnormX}
   \tnorm{x}_\nX:=\sup_{n\ge 1}\normX{P_{n} x}
\eeqal
satisfies the properties (\ref{eq:conditionX1})-(\ref{eq:conditionX3}).
\end{corollary}
\begin{proof}
It is well-known that $\tnorm\cdot_\nX$
is an equivalent norm to $\normX\cdot$, see \cite{Fabian11}. 
We have the following bound for the projections:
\beqa
  \tnorm{P_nx}_\nX &=& \sup_{1\le m<\infty}\normX{P_{m}P_n x} = \sup_{1\le m\le n}\normX{P_m x}\le \sup_{1\le m< \infty}\normX{P_{m} x} =\tnorm{x}_\nX.
\eeqa
The previous lemma now gives the assertion.    
\end{proof}
\subsection{Sum of a linear and nonlinear operator}

\begin{theorem}\label{thm:3}
We assume $\nX$ to be a separable, reflexive, real Banach space with strictly convex dual and a Schauder basis. 
Furthermore, let $\nY$ be a separable, reflexive, real Banach space.
We consider operators of the form $A= L+N:\nX\to\nY'$ with an isomorphism $L\in\mathcal L(\nX,\nY')$ and 
$N:\nX\to \nY'$ sequentially continuous for the weak topologies and it exists $0<\alpha<1$ and $c>0$ s.t.
\begin{equation}\label{eq:growth_th:3}
    \norm{N(x)}_{\nY'} \le c(1+\normX{x}^{\alpha}).
\end{equation}
Then $A$ is surjective.
\end{theorem}
\begin{proof}
Due to the assumptions on $L$ and $N$, $A$ is weakly continuous. It remains to show that
$A$ is mapped coercive (Def. \ref{def:3}), so that Theorem \ref{thm:1} gives surjectivity.
Let $\Psi:\nX\to\nX'$ be the duality map  according to Corollary \ref{lemma:dualitymap}. 
We define $S:=(L^{-1})^*\in \mathcal L(\nX',\nY)$ and obtain with (\ref{eq:psicoercive})
$$
\dpY{Lx}{S \Psi(x)} = \dpX{\Psi(x)}{S^*Lx} = \dpX{\Psi(x)}{x} \ge c_1\normX{x}^2.
$$
This gives
$$
\scp{A(x)}{S \Psi(x)} \ge c_1\normX{x}^2 - \norm{N(x)}_{\nY'}\norm{S}_{\nY\leftarrow\nX'}C_\Psi\normX{x}.
$$
Now we get with (\ref{eq:growth_th:3}) and $0<\alpha<1$:
$$
 \frac{\scp{A(x)}{S \Psi(x)}}{\norm{x}_\nX} \ge c_1\norm{x}_\nX - c'(1+\norm{x}_\nX^{\alpha}),
$$
which ensures the mapped coercivity (\ref{eq:H2XXX}). The assertion follows from Theorem~\ref{thm:2}.
\end{proof}

\section{Applications}\label{sec:5}

In this section we give several applications of Theorems \ref{thm:1} and
\ref{thm:3} for nonlinear partial differential equations in bounded Lipschitz domains $\Omega\subset\nR^d$, in $d\ge2$ dimensions. 

Applications in the context of partial differential equations suggest the setting of a
compact embedding $\nX\!\subset\!\subset\!\nZ$ with another Banach space $\nZ$.  
In this case, the weak convergence $x_n\rightharpoonup x$ in $\nX$ implies strong convergence $x_n\to x$ in $\nZ$.
This fact may help to verify the weak continuity (\ref{eq:H1}), see
examples in Sections \ref{sec:5-1} and \ref{sec:5-3}.

\subsection{Semi-linear Poisson equation in $W^{1,p}(\Omega)$}\label{sec:5-1}

Let $d\in\mathbb N$ and $\Omega\subset \nR^d$ 
be a convex bounded domain.
We consider:
\beqal\label{eq:laplace1new}
  -\Delta u  &=& f(u)\quad\mbox{in }\Omega,\qquad u = 0\quad\mbox{on }\partial\Omega.
\eeqal
Let $1<p<\infty$, $\frac1p+\frac1q=1$, $\nX=W^{1,p}_0(\Omega)$, $\nY=W^{1,q}_0(\Omega)$, such that 
$\nX'=W^{-1,q}(\Omega)$ and $\nY'=W^{-1,p}(\Omega)$. 

We will use Theorem~\ref{thm:3} to show existence of solutions.
It is well-known that $\nX$ has a strictly convex dual.
We have $L:\nX\to\nY'$ the Laplace operator with homogeneous Dirichlet boundary conditions, which is 
an isomorphisme, see Corollary~1 in \cite{Fromm1993}.
We make the assumption that $f\in C(\nR)$, with a growth condition of the type
\beqal\label{eq:f-growth}
\vert f(t)\vert \le C_H\left( 1+\vert t\vert^{\alpha}\right) \quad\forall t\in \nR,
\eeqal
with $\alpha>0$.
In what follows we use the Sobolev embedding
\beqal\label{eq:embeddings}
   \nX=W^{1,p}_0(\Omega)\subset W^{1,p}(\Omega)\hookrightarrow L^{p^{*}}(\Omega),
\eeqal
with $\frac1{p^*}=\frac1p-\frac1d$ for $p<d$, $p^*\in (1,\infty)$, if $p=d$, and $p^*=+\infty$, if $p>d$.
\begin{lemma}\label{lemma:boundflp}
Let $f:\nR\to\nR$ be bounded as (\ref{eq:f-growth}).
Then it follows that
$\norm{f(u)}_{\nY'}\le c+b\normX{u}^\alpha$ for all $u\in\nX$.
\end{lemma}
\begin{proof}
By Sobolev embedding $\nY=W^{1,q}(\Omega)\hookrightarrow L^{q^*}(\Omega)$ we get
$
  \norm{f(u)}_{\nY'} \le c\norm{f(u)}_{L^{r}}.
$
The assumed bound of $f$ implies (with different constants $b,c$)
\beqa
  \norm{f(u)}_{L^{r}} &\le& c+b\norm{u}_{L^{(1-\epsilon) r}}^\alpha.
\eeqa
Due to $1/r=1/p+1/d\ge 1/p-1/d=1/p^*$ we deduce $r\le p^*$ so that the assertion follows
by a second use of Sobolev embedding (\ref{eq:embeddings}).
\end{proof}
\begin{proposition}\label{lemma:xyz}
Under the conditions $f\in C(\nR)$ with growth condition (\ref{eq:f-growth}), $1<p<+\infty$, and $0<\alpha<1$ there exists a solution $u\in W^{1,p}_0(\Omega)$ of (\ref{eq:laplace1new}).
\end{proposition}
\begin{proof}
We will check the conditions to apply Theorem~\ref{thm:3}. The operator $A$ is the sum of $-\Delta$ and a nonlinear operator $N:\nX\to\nY'$, given by $N(u):=-f\circ u$.\\
(a) According to Fromm \cite{Fromm1993} and using the convexity of $\Omega$, $-\Delta\in L(\nX,\nY')$ is an isomorphism.\\
(b) Growth property (\ref{eq:growth_th:3}) of $N$: First, we obtain with  $\frac1{q^*}=\frac1q-\frac1d$ and $\frac1r +  \frac1{q^*} =1$ and elementary calculus that $\alpha r < r< p^*$. This allows us to derive with (\ref{eq:f-growth}) and (\ref{eq:embeddings}) for $\phi\in\nY$:
\beqa
\vert \scp{f(u)}{\phi}\vert&\le&
\norm{f(u)}_{L^{r}}\norm{\phi}_{L^{q^*}}
\le C_H\left( 1 + \norm{u}_{L^{\alpha r}}^{\alpha}\right) \norm{\nabla\phi}_{L^q}.
\eeqa
Hence, $f(u)\in \nY'$ and (\ref{eq:growth_th:3}) holds for $N(u)=f\circ u$.\\
(c) The continuity of the nonlinear part is obtained as follows: 
Let $u_n\rightharpoonup u$ in $\nX$. With the compact embedding $W^{1,p}(\Omega)\subset\subset L^{r}(\Omega)$ (because of $r<p^*$) we get $u_n\to u$ in $L^{r}(\Omega)$. By continuity of $f$ we get $f(u_n)\to f(u)$ a.e. in $\Omega$.
We  will show that this sequence is uniformly integrable over $\Omega$, i.e. $\forall\epsilon>0$ $\exists\delta>0$ s.t.
\beqal\label{eq:vitali}
A\subseteq\Omega\mbox{ measurable},\ \vert A\vert\le\delta \quad&\Longrightarrow&\quad \norm{f(u_n)}_{L^r(A)}^r\le \epsilon.
\eeqal
For $R>0$ we use the notation $A_R:=\{x\in A: \vert u_n(x)\vert \ge R\}$ and $M_R:=\!\!\max\limits_{-R\le t\le R}\vert f(t)\vert^r$. Then with (\ref{eq:f-growth})
\beqa
\norm{f(u_n)}_{L^r(A)}^r 
&=& \norm{f(u_n)}_{L^r(A_R)}^r + 
\norm{f(u_n)}_{L^r( A\setminus A_R)}^r\\
&\le& C_1 \int_{A_R} (1 + \vert u_n(x)\vert^{\alpha r}) \,dx
+ \vert A\vert M_R\\
&\le& C_1\int_{A_R} R^{\alpha r-p^*}\vert u_n(x)\vert^{p^*} \,dx
+ \vert A\vert \left(C_1 + M_R\right)\\
&\le& C_1 R^{\alpha r-p^*}\norm{u_n}_{L^{p^*}(\Omega)}+ \vert A\vert \left(C_1 + M_R\right).
\eeqa
For $ R\ge(\epsilon/(2C_1\norm{u_n}_{L^{p^*}(\Omega)}))^{p^*-\alpha r}$, and $\delta:=\epsilon/(2\left(C_1\! +\! M_R\right))$ we obtain (\ref{eq:vitali}).
By Vitali's integrability theorem \cite{BoccardoCroce2014} now gives $f(u_n)\to f(u)$ in $L^r(\Omega)$, which implies
$f(u_n) \stackrel{*}{\rightharpoonup} f(u)$ in $\nY'$.
 \end{proof}
\begin{remark}
Clearly, the existence of solutions cannot be extended to growth conditions (\ref{eq:f-growth}) with $\alpha\ge1$, since there is no solution to $-\Delta u = \lambda u$, if $\lambda\not\in \sigma(-\Delta)$. However, the nonlinear part $N(u):=f\circ u$ remains sequentially continuous for continuous $f$ with growth condition with exponent $\alpha<p^*/r=(d+p)/(d-p)$.
\end{remark}
In the Hilbert-space case $p=q=2$, we may slightly relax the growth condition:
Let $\chi_+$ and $\chi_-$ be the characteristic functions for the sets $\nR_+$ and $\nR_-$, respectively; $f_+=\max\{f,0\}$, $f_-=\min\{f,0\}$.
Let $f_0:\nR\to\nR$ defined by $f_0:=f_+\chi_++f_-\chi_-$. 
\begin{proposition}\label{prop:p2}
Let $p=q=2$, $f\in C(\nR)$  fulfills the growth condition (\ref{eq:f-growth}) with $\alpha<(d+p)/(d-p)$  for $d>2$, $1\le \alpha<\infty$ for $d\le2$, s.t. $f_0$ fulfills a similar growth barrier (\ref{eq:f-growth}) with exponent $0\le\beta<1$.
Then there exists a solution $u\in\nH:=H^1_0(\Omega)$ of (\ref{eq:laplace1new}).
\end{proposition} 
\begin{proof}
(a) We argue as in the proof of Proposition~\ref{prop:p2} with $p=q=2$.
In order to ensure  that $\vert \scp{f(u)}{\phi}\vert$ is finite and that $u_n\to u$ implies that $f(u_n)$ converges weak-$*$ to $f(u)$ it was sufficient to have $\alpha r<p^*$. This is indeed the case for 
$$
\alpha \le p^*\frac{p^*-1}{p^*} = \frac{d+2}{d-2}.
$$\\
(b) For showing the mapped coercivity (Def. \ref{def:1}) with $\Phi={\rm id}$ we use the property
 $f_0(t)=f(t)$,  if $\sgn{f(t)}=\sgn t$, and $f_0(t)=0$ otherwise.
With the Poincare inequality we get the upper bound
\beqa
  \scp{f(u)}u &=& \int_{\Omega}f(u)u\ \le\  \int_{\Omega}f_0(u)u\\
  &\le& c(\norm{u}_{L^1(\Omega)}+\norm{|u|^{1+\beta}}_{L^1(\Omega)})\\
  &\le& c'(\norm{\nabla u}_{L^2(\Omega)}+\norm{\nabla u}_{L^2(\Omega)}^{1+\beta}).
\eeqa
We use the same notation $N$ as in the proof of Proposition \ref{lemma:xyz}
and we use the fact that in the Hilbert space case we have $S\circ\Psi={\rm id}$. Therefore,
\beqa
  \scp{A(u)}u &=& \norm{\nabla u}_{L^2(\Omega)}^2 -(f(u),u)\ \ge\  \norm{\nabla u}_{L^2(\Omega)}^2 -(f_0(u),u).
\eeqa
The negative contribution $-(f_0(u),u)$ is bounded by (\ref{eq:f-growth}) as: $(1/p^*=1/2-1/d)$
\beqa
 \vert (f_0(u),u)\vert&\le& 
  c(1+\norm{|u|^{\beta}}_{L^1(\Omega)})\norm{u}_{L^1(\Omega)}\\
  &\le& c'(\norm{\nabla u}_{L^2(\Omega)}+\norm{\nabla u}_{L^2(\Omega)}^{1+\beta}),
\eeqa
where we used $1+\beta<2<p^*$ and the embedding $\nH\hookrightarrow L^{p^*}(\Omega)$.
This implies (\ref{eq:H2XX}):
\beqa
   \frac{\scp{A(u)}u }{\norm{\nabla u}_{L^2(\Omega)}} &\ge& \norm{\nabla u}_{L^2(\Omega)}-c'(1+\norm{\nabla u}_{L^2(\Omega)}^{\beta}).
\eeqa
The weak convergence  is obtained as before in the proof of Proposition \ref{lemma:xyz}.
\end{proof}

\subsection{Mixed formulation of semi-linear Poisson problem}\label{sec:5-2}

Here we give an example of linearly mapped coercivity (Def. \ref{def:1}) with a non-trivial linear $\Phi$. We may formulate (\ref{eq:laplace1new}) in the primal-mixed formulation:
\beqa
   q-\nabla u &=& 0\quad\mbox{in }\Omega,\\
   -\mbox{div}\, q &=& f(u)\quad\mbox{in }\Omega,
\eeqa
and homogeneous Dirichlet conditions.
Now, the spaces are $\nX=\nY=L^2(\Omega)\times H^1_0(\Omega)$ and the corresponding operator $A:\nX\to\nX'$ is given by
\beqa
   \dpY{A(q,u)}{(p,v)} &:=& (q,p)-(\nabla u,p) + (q,\nabla v)-(f(u),v).
\eeqa
The operator $\Phi\in  L(\nX,\nX)$ is defined as $\Phi(q,u):=(q-\nabla u,u)$.
We arrive at
\beqa
  \dpY{A(q,u)}{\Phi(q,u)} &=& \norm{q}_{L^2(\Omega)}^2-(q,\nabla u)-(\nabla u,q-\nabla u)+(q,\nabla u) -(f(u),u)\\
  &\ge&\frac12 \norm{q}_{L^2(\Omega)}^2 +\frac12\norm{\nabla u}_{L^2(\Omega)}^2-(f(u),u).
\eeqa
By the assumptions and technique as in Proposition \ref{prop:p2} for bounding $(f(u),u)$ we obtain (\ref{eq:H2XX}).
Moreover, weak continuity is ensured as shown in Section \ref{sec:5-1}.

\subsection{Navier-Stokes with Dirichlet data}\label{sec:5-3}

\quad The Navier-Stokes system in a bounded Lipschitz domain $\Omega\subset\nR^d$, $1\le d\le 3$,
with homogeneous Dirichlet data reads
\beqal\label{eq:nse1}
   (u\cdot\nabla) u -\nu\Delta u+\nabla p &=& f\quad\mbox{in }\Omega,\\
   \mbox{div}\,u &=& 0\quad\mbox{in }\Omega,\label{eq:nse2}\\
   u&=& u_0\quad\mbox{on }\partial\Omega.\label{eq:nse3}
\eeqal
Here, the quantity $u$ is vector-valued. 
Accordingly, the $L^2$-inner product used in the following is the vector-valued inner product. We now present a short proof of existence of weak solutions without smnallness assumption on the data $u_0,f$.
We firstly treat the case of homogeneous Dirichlet data, i.e. $u_0=0$.
\begin{corollary}\label{cor:nse1}
In a bounded Lipschitz domain $\Omega\subset\nR^d$, $1\le d\le3$,
the Navier-Stokes problem (\ref{eq:nse1})-(\ref{eq:nse3}) with $u_0=0$ has for each $f\in H^{-1}(\Omega)^d$
a weak solution 
$(u,p)\in H_0^1(\Omega)^d\times$ $L^2(\Omega)$.
\end{corollary}
\begin{proof}
 The variational space for the velocity is given by the divergence free $H^1$-func\-tions:
\beqa
   \mathcal H_0:=\{u\in H_0^1(\Omega)^d : (\mbox{div}\,u,q)=0\ \forall q\in L^2(\Omega)\}.
\eeqa
The corresponding operator for the velocities reads for $u\in\nH_0$:
\beqa
  A(u) &:=& (u\cdot\nabla)u - \nu\Delta u.
\eeqa
Existence of $u\in \mathcal H_0$ follows from Lemma \ref{lemma:1} with $\nX=\nY=\nH$:
Weak continuity is verified similar to the previous example by showing convergence of the nonlinear term
$
   \lim_{n\to\infty}((u\cdot\nabla)u-(u_n\cdot\nabla)u_n,\phi)= 0
$
for $u_n\to u$ in $L^2$ and $u_n\rightharpoonup u$ in $\nH$. To this end, one considers firstly
$\phi\in\mathcal D(\Omega)$ and secondly $\phi\in\mathcal H_0$, see Lemma 1.5 in \cite{Temam1984} for details.
Coercivity (\ref{eq:standardcoercivity}) is ensured, because of $((u\cdot\nabla)u,u)=0$.

The existence of the
pressure $p$ follows by the inf-sup property of the divergence operator, see Temam \cite{Temam1984} for details.
\end{proof}
We like to refer to a different proof of Navier-Stokes solutions in $W^{1,p}(\Omega)$ with $3<p$ $<\infty$ by Benjamaa et al. \cite{Benjemaa2017}.
Now, we treat the case of inhomogeneous Dirichlet data.
\begin{corollary}\label{cor:nse2}
In a bounded Lipschitz domains $\Omega\subset\nR^d$, $1\le d\le3$,
the Navier-Stokes problem (\ref{eq:nse1})-(\ref{eq:nse3}) has for each $u_0\in H^s(\partial\Omega)^d$,
$s\in(\frac12,1)$,
and each $f\in H^{-1}(\Omega)^d$
a weak solution $(u,p)\in H^1(\Omega)^d\times L^2(\Omega)$.
\end{corollary}
\begin{proof}
For the Navier-Stokes system (\ref{eq:nse1})-(\ref{eq:nse3}) with inhomogeneous Di\-richlet data the Hilbert space is here given by 
\beqa
   \mathcal H:=\{u\in H^1(\Omega)^m : (\mbox{div}\,u,q)=0\ \forall q\in L^2(\Omega)\}.
\eeqa
In order to formulate the corresponding variational formulation, one may choose  
an arbitrary extension $u_{ext}\in \mathcal H$ of the Dirichlet data, i.e. $u_{ext}\vert_{\partial\Omega}=u_0$.
Making the ansatz $u+u_{ext}$ and seeking $u\in \mathcal H_0$ as the solution of
the corresponding Operator $A:\mathcal H_0\times\mathcal H_0\to (\mathcal H_0\times\mathcal H_0)'$:
\beqa
  A(u) &:=& ((u_{ext}+u)\cdot\nabla) (u_{ext}+u) - \nu\Delta (u_{ext}+u).
\eeqa
Coercivity can be shown by the same techniques as in the proof of Corollary \ref{cor:nse1}. 
Diagonal testing yields
\beqa
   \dpY{A(u)}{u} &=&((u\cdot\nabla)u_{ext},u) + ((u_{ext}\cdot\nabla)u_{ext},u) +\nu\norm{\nabla u}_{L^2}^2 
   + \nu(\nabla u_{ext},\nabla u).
\eeqa
For (\ref{eq:H2XX}), the linear terms are  not critical. For the nonlinear term we use the result of
Neustupa \cite{Neustupa2022}: for $u_0\in$ $H^s(\partial\Omega)^d$, $s\in(\frac12,1)$,
and arbitrary $\epsilon>0$ there can be found a particular extension $u_{ext}^*\in \mathcal H_0$ of $u_0$
such that
\beqa
  \vert ((u\cdot\nabla)u_{ext}^*,u) \vert&\le& \epsilon\norm{\nabla u}_{L^2(\Omega)}^2\quad\forall u\in \mathcal H_0.
\eeqa
As a consequence, there is an extension $u_{ext}^*$ s.t. we obtain coercivity (\ref{eq:standardcoercivity}). 
This ensures existence of solutions in the non-homogeneous case.
\end{proof}

\normalsize

\bibliographystyle{plain} 
\bibliography{mybib}

\end{document}